\documentclass[psamsfonts, 12pt]{amsart}
\usepackage{graphicx, psfrag}

 \newcommand{\R}{\Bbb R}
 \newcommand{\Z}{\Bbb Z}
\newcommand{\ep}{\varepsilon}
\DeclareMathOperator{\Int}{Int} \DeclareMathOperator{\Fix}{Fix}

\newcommand{\nw}{\Omega}

\theoremstyle{plain} \newtheorem{thm}{Theorem}
\newtheorem{cor}[thm]{Corollary} \newtheorem{prop}[thm]{Proposition}
\newtheorem{lemma}[thm]{Lemma}

\theoremstyle{definition} \newtheorem{defn}[thm]{Definition}

\newtheorem{ex}[thm]{Example} 

\theoremstyle{remark}

\begin{document}

\title{Bounded homeomorphisms of the open annulus}

\author{David Richeson} \author{Jim Wiseman} \address{Dickinson
College\\Carlisle, PA 17013} \email{richesod@dickinson.edu}
\address{Swarthmore College\\Swarthmore, PA 19081}
\thanks{The second author was partially supported by the Swarthmore College
Research Fund.}

\email{jwisema1@swarthmore.edu}
% 37E30

\date{\today}
\begin{abstract} We prove a generalization of the Poincar\'e-Birkhoff theorem
for the open annulus showing that if a homeomorphism satisfies a certain twist
condition and the nonwandering set is connected, then there is a fixed point. 
Our main focus is the study of bounded homeomorphisms of the open annulus.  We
prove a fixed point theorem for bounded homeomorphisms and study the special
case of those homeomorphisms possessing at most one fixed point.  Lastly we use
the existence of rational rotation numbers to prove the existence of periodic
orbits.
\end{abstract}

\maketitle
\section{Introduction} A homeomorphism $f:X\to X$ is said to be \emph{bounded}
if there is a compact set which intersects the forward orbit of every point. 
Since every homeomorphism on a compact space is bounded, bounded homeomorphisms
are interesting only on noncompact spaces.  As we will see, if $f$ is bounded
then there is a forward invariant compact set which intersects the forward
orbit of every point.  Thus, a bounded map on a noncompact space behaves in
many ways like a map on a compact space.  In particular, many results that are
true for maps on compact spaces are also true for bounded maps on noncompact
spaces (e.g., the Lefschetz fixed point theorem).  

In this paper we study primarily the dynamics of bounded homeomorphisms of the
open annulus.  Intuitively we may view these homeomorphisms as those having
repelling boundary circles.  In fact, we will see that the orbit of every point
intersects an essential, closed, forward invariant annulus.  Thus, roughly
speaking, many of the results for homeomorphisms of the closed annulus also
hold for bounded homeomorphisms of the open annulus.   Conversely, many of the
results that hold for bounded homeomorphisms of the open annulus also hold for
homeomorphisms of the closed annulus; one may enlarge the closed annulus to an
open annulus and extend the homeomorphism to a bounded homeomorphism of this
open annulus. 

The most celebrated result for the closed annulus is the Poincar\'e-Birkhoff
theorem (also called Poincar\'e's last geometric theorem), which states than
any area preserving homeomorphism which twists the boundary components in
opposite directions has at least two fixed points.  In \cite{F2} Franks gives a
topological generalization for the open annulus; he proves that if every point
in an open annulus is nonwandering and $f$ satisfies a twist condition, then
there is a fixed point of positive index.  We prove a further generalization
showing that if $f$ satisfies a twist condition and the nonwandering set is
connected then $f$ has a fixed point.  Recall that for a map $f:X\to X$, a
point $x\in X$ is {\em nonwandering} if for every open set $U$ containing $x$
there exists $n>0$ such that $f^n(U)\cap U\ne\emptyset$.  The collection of
nonwandering points is the {\em nonwandering set}, denoted $\nw(f)$. 

The paper is divided as follows.  In Section \ref{sec:bdded} we present general
properties of bounded homeomorphisms of the annulus.  In Section
\ref{sec:pbthm} we prove a generalization of the Poincar\'e-Birkhoff-Franks
theorem for the open annulus.  This section applies to homeomorphisms of the
open annulus that need not be bounded.  It can be read independently of the
rest of the paper and may be of more general interest.  In Section
\ref{sec:bddfixed} we use this theorem to prove a fixed point theorem for
bounded homeomorphisms of the open annulus.  It is interesting to note that a 
bounded homeomorphism of a noncompact space can never preserve Lebesgue
measure.  Thus, we prove a fixed point theorem for a family of maps far from
satisfying the hypotheses of the Poincar\'e-Birkhoff theorem.  Also, we study
the special case of those bounded homeomorphisms having at most one fixed
point.   Lastly, in Section \ref{sec:periodic} we apply the theorem to those
bounded homeomorphisms having a point with a rational rotation number and prove
the existence of a periodic point with that same rotation number.

In this paper we will let $A$ denote the annulus $(\R/\Z)\times I$, where
$I=[0,1]$ if $A$ is the closed annulus, and $I=(0,1)$ if $A$ is the open
annulus.  $\tilde A=\R\times I$ will denote the universal cover of the annulus
$A$ with  $\pi:\tilde A\to A$ being the covering projection.  We view $\tilde
A$ as a subset of $\R^2$, thus when we subtract two elements in $\tilde A$ we
obtain a vector in $\R^2$.  The projection onto the first coordinate
$\R^2\to\R$ is given by $(x,y)_1= x$. For any set $U\subset \tilde A$, let
$U+k$ denote the set $\{(x+k,y)\in\tilde A:(x,y)\in U\}$.  

If $f:A\to A$ is a homeomorphism then there is a {\em lift}, $\tilde f:\tilde
A\to\tilde A$ satisfying $\pi\circ\tilde f=f\circ\pi$.  Notice that $\tilde g$
is another lift of $f$ iff $\tilde g(x,y)=\tilde f(x,y)+(k,0)$ for some integer
$k$.  For any $y\in\tilde A$ define $\rho(y,\tilde f)$ to be
$\displaystyle\lim_{n\to\infty}(1/n)(\tilde f^n(y)-y)_1$ (if this limit
exists).  If $\tilde g$ is another lift then $\rho(y,\tilde g)=\rho(y,\tilde
f)+k$ for some integer $k$.  Thus we may define the {\em rotation number} of
$x=\pi(y)\in A$ to be $\rho(x)=\rho(y,\tilde f)$ (mod 1) if this limit exists. 
So defined, $\rho(x)$ is independent of the choice of $y$ and $\tilde f$. 
Unlike the case of homeomorphisms of the circle, for homeomorphisms of the
annulus different points may have different rotation numbers, and it may happen
that the rotation number for a point does not exist.

\section{Bounded homeomorphisms of the annulus}
\label{sec:bdded}

ln \cite{RW} the authors introduced the following definitions.

\begin{defn} A compact set $W$ is a {\em window} for a dynamical system on $X$
if the forward  orbit of every point $x\in X$ intersects $W$.  If a dynamical
system  has a window then we will say that it is {\em bounded}.
\end{defn}

We showed that we can characterize bounded dynamical systems in many ways.  The
following theorem summarizes some results from \cite {RW}.

\begin{thm}\label{thm:equivalences} Let $f:X\to X$ be a continuous map on a
locally compact space $X$.  Then the following are equivalent.
\begin{enumerate}
\item
$f$ is bounded.
\item There is a forward invariant window.
\item Given any compact set $S\subset X$ there is a window $W\subset X$ 
containing $S$ such that $f(W)\subset\Int W$.
\item There is a compact set $W\subset X$ with the property that 
$\emptyset\ne\omega(x)\subset W$ for all $x\in X$.
\item
$f$ has a compact global attractor $\Lambda$ (i.e., $\Lambda$ is an  attractor
with the property that for every $x\in X$, $\omega(x)$ is nonempty and
contained in
$\Lambda$).
\end{enumerate}
\end{thm}

Because every bounded map has a compact global attractor it is impossible for
it to preserve Lebesgue measure on a noncompact space.  Thus we have the
following corollary.

\begin{cor}\label{cor:areapres} Suppose $f:X\to X$ is an area preserving map of
a noncompact space $X$.  Then $f$ is not bounded.  In particular, if $S\subset
X$ is any compact set, then there exists a point $x\in X$ such that the forward
orbit of $x$ does not intersect $S$.
\end{cor}

\begin{ex}\label{ex:billiards} Consider a convex billiards table.  Is it
possible to find a  trajectory with the property that the angle the ball makes
with the  wall is always smaller than some arbitrarily chosen $\ep$? We  see
that the answer is yes.

Let $f:S^1\times(0,\pi)\to S^1\times(0,\pi)$ be the billiards map 
corresponding to the given table.  It is well known that $f$ is an  area
preserving homeomorphism homotopic to the identity.  By  Corollary
\ref{cor:areapres} $f$ is not bounded.  In particular,  there exists a point
$(x,\theta)$ whose forward orbit does not  intersect the closed annulus
$S^1\times[\ep,\pi-\ep]$.

Thus, for any $\ep>0$, there exists a trajectory 
$(x_0,\theta_0), (x_1,\theta_1), (x_2,\theta_2),...$ such that either 
$\theta_k<\ep$ for all $k\ge 0$ or $\pi-\theta_k<\ep$ for all $k\ge 0$.
\end{ex}

\begin{ex} Suppose there is a convex billiards table with bumpers in the middle
of the table (see Figure \ref{fig:billiards}).  Is it possible to find a
trajectory of the billiards ball that never strikes a bumper?

Assume that the bumpers are a finite collection of compact sets not touching
the wall of the billiards table.  Consider the billiards map for the table with
no bumpers, $f:S^1\times(0,\pi)\to S^1\times(0,\pi)$.  Let $W\subset
S^1\times(0,\pi)$ be the set of points $\{(x,\theta)\}$ with the property that
a ball at position $x$ with trajectory angle $\theta$ will strike a bumper
before striking the wall again.  Clearly $W$ is a compact set.  Thus, we
rephrase the question:  Is it possible to find an orbit of $f$ that never
intersects $W$?  By the discussion in Example \ref{ex:billiards} it is clear
that such a trajectory does exist.  Thus, given any compact set of bumpers,
there is always a trajectory that avoids the bumpers.
\begin{figure}
	\centering
	\includegraphics{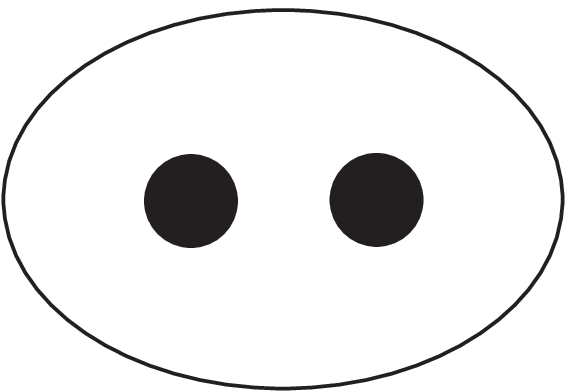}
	\includegraphics{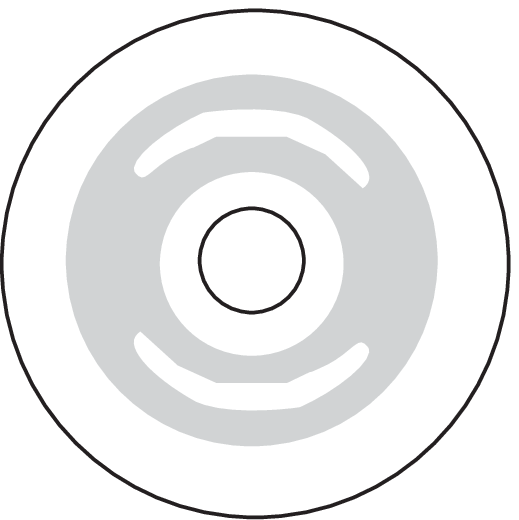}
	\caption{The billiards table with two bumpers and the corresponding
configuration space}
	\label{fig:billiards}
\end{figure}
\end{ex}

\begin{prop}\label{prop:connected} Suppose $f:A\to A$ is a bounded
homeomorphism of the open annulus  with a compact global attractor
$\Lambda\subset A$.  Then the  following are true.
\begin{enumerate}
\item The inclusion $i:\Lambda\to A$ induces an isomorphism on \v{C}ech 
cohomology, $i^{*}:\Check{H}^{*}(A)\to\Check{H}^{*}(\Lambda)$.
\item
$\Lambda$ is connected.
\item
$\Lambda$ separates the two boundaries of $A$.
\end{enumerate}
\end{prop}
\begin{proof} Let $f:A\to A$ be a bounded homeomorphism of the open annulus $A$ 
with compact global attractor $\Lambda$.  By Theorem 
\ref{thm:equivalences} there exists a window $W$ such that 
$\Lambda\subset f(W)\subset\Int W$.  Let $\ep>0$ be small  enough such that
$\Lambda\subset  A_\ep=[\ep,1-\ep]$.  For each $x\in  A_\ep$ there exists
$n_x>0$ such that $f^{n_x}(x)\subset\Int  W$.  There exists an open set $U_x$
containing $x$ such that 
$f^{n_x}(U_x)\subset\Int W$.  The collection $\{U_x\}$ is an open  cover of
$A_\ep$, thus there exists a finite subcover, 
$\{U_{x_1},\dots,U_{x_m}\}$.  Let $N=\max\{n_{x_1},\dots,n_{x_m}\}$.  It
follows that $f^i(A_\ep)\subset\Int W$ for all $i\ge N$.

Notice that $f^{N}(A_\ep)$ separates the two boundaries of 
$A$ and $f^N$ induces an isomorphism on cohomology.  Also, 
$U=\Int(A_\ep), f^{N}(U), f^{2N}(U),\ldots$ is a nested  sequence of open sets
with $\Lambda=\bigcap_{k=0}^{\infty}f^{kN}(U)$.  Consequently, the inclusion
$i:\Lambda\to A$ induces an isomorphism 
$i^{*}:\Check{H}^{*}(A)\to\Check{H}^{*}(\Lambda)$ and $\Lambda$  separates the
two boundaries of $A$.  Moreover, since $\Lambda$ is  the intersection of a
nested collection of connected open sets, 
$\Lambda$ is itself connected.
\end{proof}

Next we prove a key result that states that all of the interesting dynamics
occurs inside a closed annulus.  This result is very useful.  It validates our
intuition that a bounded homeomorphism on the open annulus behaves like a
homeomorphism on the closed annulus.

\begin{prop}\label{prop:closedannulus} If $f:A\to A$ is a bounded homeomorphism
of an open annulus, then  there exists a closed annulus $A_0\subset A$ whose
boundaries are  smooth essential curves such that $f(A_0)\subset\Int A_0$. 
Moreover, 
$A_0$ can be chosen so that  boundary is as close to $\Lambda$ or as  close to
the boundary of $A$ as desired.
\end{prop}
\begin{proof} Let $f:A\to A$ be a bounded homeomorphism of the open annulus 
$A=S^1\times (0,1)$.  By Theorem \ref{thm:equivalences} there exists  a compact
global attractor $\Lambda\subset A$.  Let $\ep>0$ ($\ep$ should be small enough
that $\Lambda\subset S^1\times[\ep,1-\ep]$).  We will construct a closed
annulus $A_0$ satisfying the conclusion of  the theorem with the property that 
$[\ep,1-\ep]\subset A_0$.  A similar argument can be  used to show that we can
find $A_0$ with the boundary close to 
$\Lambda$.

Let $A^*=A\cup\{*\}$ be the one point compactification of $A$.  It is  easy to
see that $(\Lambda,\{*\})$ is an attractor-repeller pair (in  the sense of
Conley \cite{C}).  Let  $\gamma:A^*\to\R$ be a  continuous Lyapunov function
satisfying $\gamma^{-1}(0)=\Lambda$, 
$\gamma^{-1}(1)=\{*\}$ and $\gamma(f(x))<\gamma(x)$ for all 
$x\not\in(\Lambda\cup\{*\})$ (see \cite{F1} for details). For the  remainder of
the proof we will restrict $\gamma$ to be a function  from $A$ to $\R$.  Let
$c\in(0,1)$ be such that 
$\gamma^{-1}(c)\cap(S^1\times[\ep/2,1-\ep/2])=\emptyset$.  Because $\gamma$ may
not be smooth the set $\gamma^{-1}(c)$ could be  quite complicated.  For any
smooth function $\lambda:A\to\R$ (which  may not be a Lyapunov function)
sufficiently $C^0$-close to $\gamma$  and any regular value for $\lambda$,
$c^\prime\in\R$, sufficiently  close to $c$, 
$\lambda^{-1}(c^\prime)\cap(S^1\times[\ep,1-\ep])=\emptyset$  and
$\lambda^{-1}(c^\prime)\cap f(\lambda^{-1}(c^\prime))=\emptyset$.  Because
$c^\prime$ is a regular value, $\lambda^{-1}(c^\prime)$ is  the disjoint union
of smoothly embedded circles in $A$. By  Proposition \ref{prop:connected},
$\Lambda$ separates the two  boundaries of $A$.  Thus there is one circle in 
$\lambda^{-1}(c^\prime)$ that separates $\Lambda$ from the inside  boundary and
another circle that separates $\Lambda$ from the outside  boundary.  The region
bounded by these two circles is a closed  annulus $A_0$ with
$[\ep,1-\ep]\subset A_0\subset A$  and $f(A_0)\subset\Int A_0$.
\end{proof}

\begin{cor}\label{cor:indexzero} If $f:A\to A$ is a bounded homeomorphism of
the open annulus homotopic to the identity, then the Lefschetz index of the
fixed point set is zero.  In particular, if
$f$ has a fixed point of nonzero index, then $f$ has at least two fixed points.
\end{cor}
\begin{proof} Suppose $f:A\to A$ is bounded.  Then there exists an essential
closed annulus 
$A_0\subset A$ containing the fixed point set with the property that 
$f(A_0)\subset\Int A_0$.  So, the fixed point set of $f$ has Lefschetz  index
zero.  Clearly, if $f$ has a fixed point of nonzero index, then 
$f$ has at least two fixed points.
\end{proof}

\section{A generalization of the Poincar\'e-Birkhoff theorem}
\label{sec:pbthm}

The classical Poincar\'e-Birkhoff Theorem states that every area preserving
homeomorphism of the closed annulus that twists the two boundary components in
opposite directions must have two fixed points (\cite{P}, \cite{B1},
\cite{B2}).  In the years since it was proved there have been new proofs and
various generalizations (see for instance \cite{BN},\cite{F2}, \cite{F4}, 
\cite{Ca}, \cite{G}, \cite{W}, \cite{AS}).  In \cite{F2} Franks generalizes
this theorem to the open annulus.  He weakens the area preserving hypothesis to
the assumption that every point is nonwandering and he weakens the twist
condition to one about positively and negatively returning disks.  The expense
of these assumptions is that the homeomorphism may have only one fixed point,
but this fixed point has positive index.

In this section we observe that we may weaken the hypotheses to the assumption
that the nonwandering set, $\nw(f)$,  is connected.  In this case the
homeomorphism must have a fixed point (now possibly of zero index). Since there
are now points that are not nonwandering we must clarify the twist condition -
we will insist that the positively and negatively returning disks intersect the
nonwandering set.

\begin{defn}\label{defn:returning} Let $f:A\to A$ be a homeomorphism of an open
or closed annulus and let $\tilde f:\tilde A\to \tilde A$ be a lift of $f$.  An
open disk $U\subset\tilde A$ is a {\em positively returning disk} if $\tilde
f(U)\cap U=\emptyset$, if $\pi(U)$ is a disk in $A$, and if there exist $n,k>0$
such that $\tilde f^n(U)\cap (U+k)\ne\emptyset$. Define the set $\nw^+(f)
=\{y\in\nw(f): y\in\pi(U)\text{ for some positively returning disk }U\}.$ 
Simlarly, define {\em negatively returning disks} (requiring $k<0$) and
$\nw^-(f)$. Notice that these definitions depend on the choice of the lift.
\end{defn}

Observe that  the nonwandering set of $f$, $\nw(f)$, is equal to the (not
necessarily disjoint) union of $\nw^+(f)$, $\nw^-(f)$, and $\pi(\nw(\tilde f))$.

We will need the following definition from \cite{F2}.

\begin{defn} Let $f:M\to M$ be a homeomorphism of a surface.  A {\em disk
chain} for $f$ is a finite collection of embedded open disks,
$U_1,...,U_n\subset M$ satisfying
\begin{enumerate}
\item
$f(U_i)\cap U_i=\emptyset$ for all $i$.
\item For all $i,j$, either $U_i=U_j$ or $U_i\cap U_j=\emptyset$.
\item For each $i<n$ there exists a positive integer $m_i$ such that
$f^{m_i}(U_i)\cap U_{i+1}\ne\emptyset$.
\end{enumerate} If $U_1=U_n$ then we say that $U_1,...,U_n$ is a {\em periodic
disk chain}.
\end{defn}

Franks proves the following generalization of a theorem of Brouwer (see also
\cite{Br}, \cite{Fa}).

\begin{thm}\label{thm:perimpliesfixed}\cite{F2} Suppose $f:\R^2\to\R^2$ is an
orientation preserving homeomorphism with isolated fixed points.  If $f$ has a
periodic disk chain, then $f$ has a fixed point of positive index. In
particular, if $f$ has a periodic point, then $f$ has a fixed point.
\end{thm}

The following two lemmas are consequences of this theorem.

\begin{lemma}\label{lem:nwfp} Suppose $f:\R^2\to\R^2$ is an orientation
preserving homeomorphism. If $\nw(f)\ne\emptyset$ then $f$ has a fixed point. 
If $\nw(f)$ consists of more than just fixed points and the fixed points are
isolated, then $f$ has a fixed point of positive index.
\end{lemma}
\begin{proof} Let $x\in\nw(f)$.  If $x$ is not a fixed point, then there exists
an open disk $U$ containing $x$ such that $f(U)\cap U=\emptyset$.  Since $x$ is
nonwandering there exists $n>1$ such that $f^n(U)\cap U\ne\emptyset$.  Thus
$U_1=U_2=U$ is a periodic disk chain.  By Theorem \ref{thm:perimpliesfixed} $f$
has a fixed point, and if the fixed points are isolated, then there is a fixed
point of positive index.
\end{proof}

Although not explicitly stated as a result, the following lemma was proved in
\cite{F2}.  

\begin{lemma}\label{lem:posandneg} Suppose $f:A\to A$ is an orientation
preserving homeomorphism of the open annulus that is homotopic to the identity,
and let $\tilde f:\tilde A\to\tilde A$ be a lift of $f$.  If there is a disk
$U\subset\tilde A$ that is both positively and negatively returning, then
$\tilde f$, and hence $f$, has a fixed point.  If the fixed points are
isolated, then there is a fixed point of positive index.
\end{lemma}
\begin{proof} Suppose $U\subset\tilde A$ is both a positively and negatively
returning disk.  So, there exist $n_{1},n_{2},k_{1},k_{2}>0$ such that $\tilde
f^{n_{1}}(V)\cap (U+k_{1})\ne\emptyset$ and $\tilde f^{n_{2}}(U)\cap
(U-k_{2})\ne\emptyset$.  As shown in \cite{F2}, $U+k_1$, $U+2k_1$,
$U+3k_1$,..., $U+k_2k_1$, $U+(k_1-1)k_2$,..., $U+2k_2$, $U+k_2$, $U$ is a
periodic disk chain.  Thus, by Theorem \ref{thm:perimpliesfixed} the
conclusions hold.
\end{proof}

We now give our main theorem of this section, a generalization of the 
Poincar\'e-Birkhoff-Franks theorem.

\begin{thm}\label{thm:connectednw} Suppose $f:A\to A$ is an orientation
preserving homeomorphism of the  open annulus that is homotopic to the
identity, and suppose the  nonwandering set of $f$, $\nw(f)$, is connected.  If
there is a lift $\tilde  f:\tilde A\to\tilde A$ possessing a positively
returning disk and a  negatively returning disk both intersecting
$\pi^{-1}(\nw(f))$, then 
$\tilde f$, and hence $f$, has a fixed point.
\end{thm}
\begin{proof} Let $f$ and $\tilde f$ be as above.  For the sake of
contradiction,  suppose $\tilde f$ has no fixed point.  Since $\tilde A$ is 
homeomorphic to $\R^2$, Lemma \ref{lem:nwfp} implies that $\nw(\tilde 
f)=\emptyset$.  By the remark following Definition
\ref{defn:returning} we know that $\nw(f)=\nw^+(f)\cup\nw^-(f)$. 
 From their definitions it is easy to see that $\nw^+(f)$ and 
$\nw^-(f)$ are open subsets of $\nw(f)$.  Since $\nw(f)$ is connected  it
follows that $\nw^+(f)\cap\nw^-(f)\ne\emptyset$.

Let $x\in\tilde A$ with $\pi(x)\in\nw^+(f)\cap\nw^-(f)$.  Then there  exists a
positively returning disk, $U_1$ and a negatively returning  disk, $U_2$, both
containing $x$.  Let $U\subset U_1\cap U_2$ be an open disk  containing $x$. 
Since $\pi(x)$ is nonwandering and since $\nw(\tilde  f)=\emptyset$ the disk
$U$ must be either positively or negatively  returning.  If it is positively
returning then $U_2$ must also be  positively returning.  Similarly, if $U$ is
negatively returning then 
$U_1$ must also be negatively returning.  Thus either $U_1$ or $U_2$  is both 
positively and negatively returning.  By Lemma 
\ref{lem:posandneg} $\tilde f$ has a fixed point.  This is a  contradiction. 
Thus $\tilde f$, and hence $f$, must have a fixed point.
\end{proof}

It is worth making a few comments about the hypotheses of Theorem
\ref{thm:connectednw}.  First of all, notice that the assumption that $\nw(f)$
is connected is stronger than we need.  If there exist positively and
negatively returning disks that intersect the same connected component of
$\pi^{-1}(\nw(f))$ then we could use the same proof to show the existence of a
fixed point.  Secondly,  in the definition of positively and negatively
returning disks we assume that $k\ne 0$ for both definitions.  One may ask if
the existence of returning disks with $k=0$ could be incorporated in Theorem
\ref{thm:connectednw}.  For instance, if there is a  homeomorphism with a
positively returning disk and a returning disk with $k=0$, is there a fixed
point?  The answer is yes; in fact, there is a fixed point even without the
positively returning disk.  If there is an open disk $U$ satisfying the
definition of the returning disks but with $k=0$ then $U_1=U_2=U$ is a periodic
disk chain and thus Theorem \ref{thm:perimpliesfixed} guarantees the existence
of a fixed point.

Unlike the Poincar\'e-Birkhoff theorem, our proof can guarantee only one fixed
point (not two).  Also, unlike in Franks' theorem, this one fixed point may
have index zero.  We have the following example showing that this may indeed
occur.  The example is based on one from Carter (\cite{Ca}).

\begin{ex}
\label{ex:onefixed} Consider the flow on $\tilde A$ shown in Figure
\ref{fig:onefixed}.  Let $\tilde f:\tilde A\to\tilde A$ be the time-one map of
this flow  and let $f:A\to A$ be the corresponding map on the open annulus.  So 
defined, $f$ is a bounded homeomorphism with only one fixed point.  By 
Corollary \ref{cor:indexzero} this fixed point must have index  zero.
\begin{figure}[ht]
	\centering
	\includegraphics[width=4in]{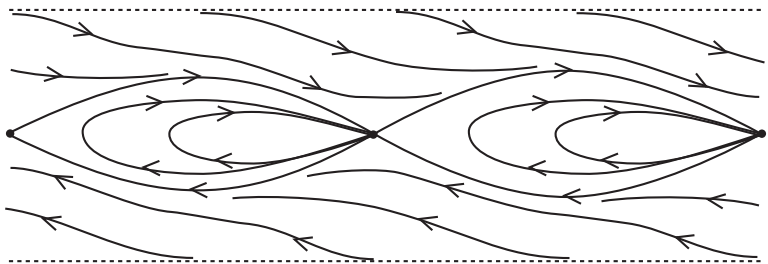}
	\caption{The lift  of a map with one fixed point of index zero}
	\label{fig:onefixed}
\end{figure}
\end{ex}

Moreover, the next example illustrates that  it is necessary for the positively
and negatively returning disks to intersect the lift of the nonwandering set.  
The positively and negatively returning disks give us reliable twist
information only if they have some recurrence.

\begin{ex}
\label{ex:posneg} Consider the time-one map, $\tilde f:\tilde A\to\tilde A$, of
the flow shown in Figure \ref{fig:posneg}.  Let $f:A\to A$ be the corresponding
map of the open annulus.  The map $f$ is a bounded homeomorphism with a
connected nonwandering set.  Moreover, $\tilde f$ possesses positively and
negatively returning disks.  Yet $f$ has no fixed point. 
\begin{figure}[ht]
	\centering
	\includegraphics[width=4in]{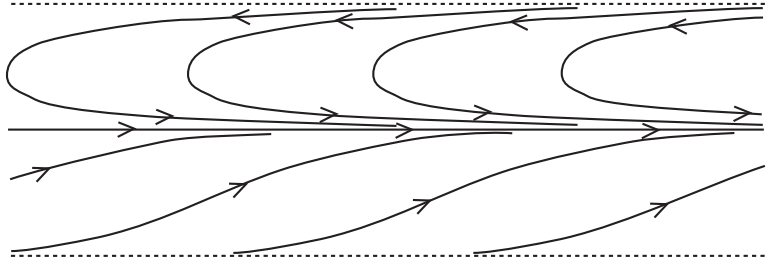}
	\caption{A map with positively and negatively returning disks and  no fixed
points}
	\label{fig:posneg}
\end{figure}
\end{ex}

In Example \ref{ex:posneg} we see that the fact that a point is in a 
negatively returning disk does not necessarily imply that the points  toward
which it tends are in negatively returning disks themselves. However, the
converse is true, as the next proposition shows.

\begin{prop} Let $f:A\to A$ be a homeomorphism of an open or closed annulus and 
let $\tilde f:\tilde A\to \tilde A$ be a lift of $f$.   Let $x\in A$.  If
$\omega(x)\cap\nw^+(f)\ne\emptyset$ then there is a positively  returning disk
containing $y\in\pi^{-1}(x)$.  If 
$\omega(x)\cap\nw^-(f)\ne\emptyset$ then there is a negatively  returning disk
containing $y\in\pi^{-1}(x)$.
\end{prop}
\begin{proof} Suppose $\omega(x)\cap\nw^+(f)\ne\emptyset$ and
$y\in\pi^{-1}(x)$.  Let $z\in\omega(x)\cap\nw^+(f)$. Then there exists a
positively  returning disk $U$ such that $z$ is in $\pi(U)$.  Also, there
exists 
$n>0$ such that $f^n(x)\in\pi(U)$.  Without loss of generality we may  assume
that $\tilde f^n(y)\in U$ (if not then translate $U$ by the  appropriate
integer amount).  Since $U$ is a positively returning  disk then so is
$V=\tilde f^{-n}(U)$.  Moreover, $V$ contains $y$.  The case for negatively
returning disks is proved similarly.
\end{proof}

In Examples \ref{ex:onefixed} and \ref{ex:posneg} we see that $\nw^+(f)$ and
$\nw^-(f)$ are disjoint sets.  Example \ref{ex:triplehorse} shows that this
need not be the case in general. Moreover, we will see that for a point $x$
with $\pi(x)\in\nw^+(f)\cap\nw^-(f)$, there may be positively and negatively
returning disks containing $x$ that are arbitrarily small.

\begin{ex}
\label{ex:triplehorse} We begin with a rectangle $N$ and create a triple
horseshoe by wrapping $N$ around the annulus twice (see Fig. 
\ref{fig:triplehorse}).  Extend $f$ to a homeomorphism on all of $A$. If
desired we may make $f$ bounded.  Choose a lift $\tilde f$ as shown in Fig.
\ref{fig:triplehorse}.

\begin{figure}[ht]
	\centering
	\psfrag{n}{$N$}
	\psfrag{fu}{$\tilde f(U)$}
	\psfrag{u}{$U$}
	\psfrag{fn}{$f(N)$}
	\psfrag{v}{$V$}
	\psfrag{n0}{$N_0$}
	\psfrag{n1}{$N_1$}
	\psfrag{n2}{$N_2$}
	\includegraphics[width=4in]{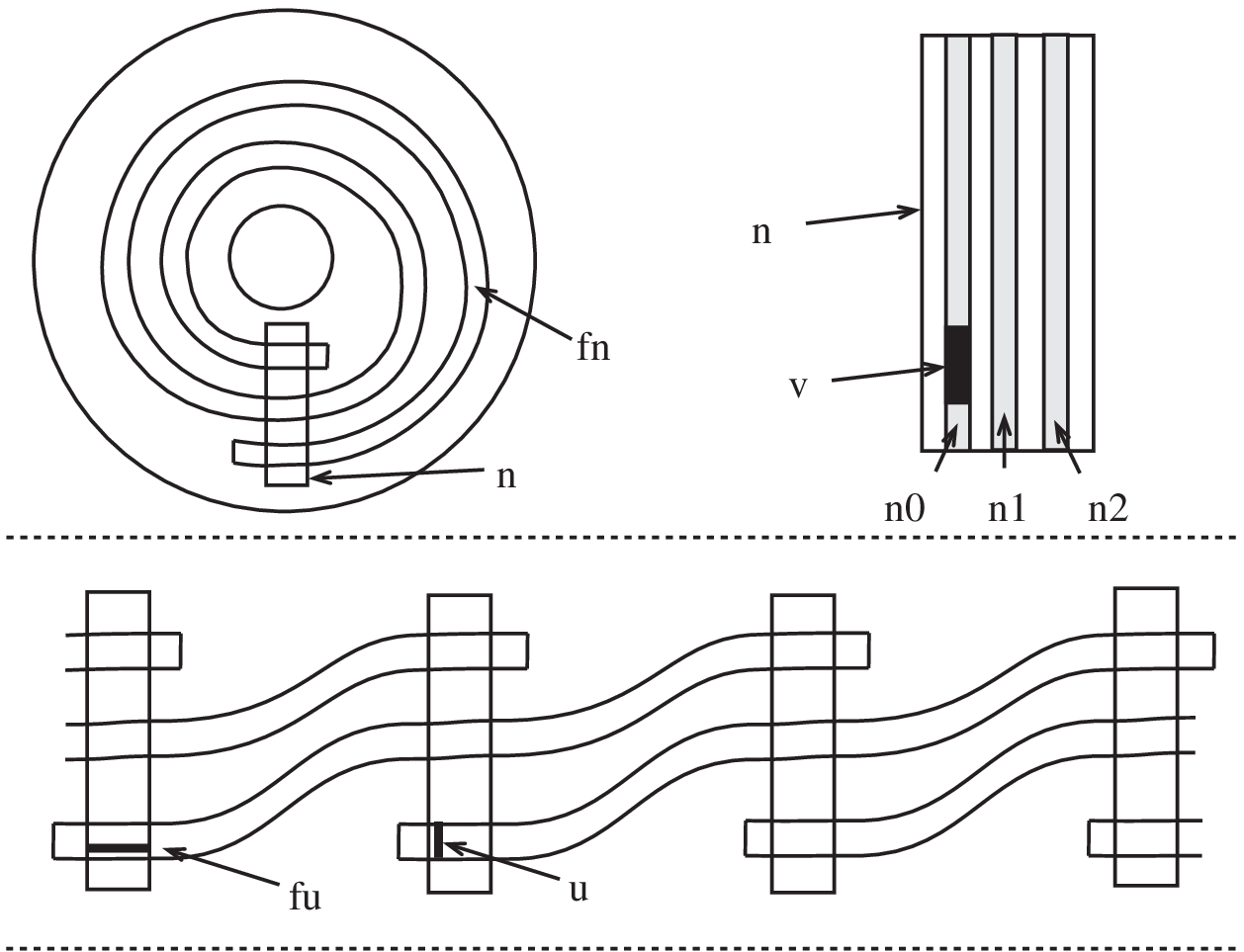}
	\caption{A map with $\nw^+(f)\cap\nw^-(f)\ne\emptyset$}
	\label{fig:triplehorse}
\end{figure} Inside this triple horseshoe is an invariant set $S$ on which $f$
is  conjugate to the full three-shift $\Sigma_3$.  In particular, let 
$N_0,N_1,N_2\subset N$ be the three components of $N\cap f^{-1}(N)$.  Then the
conjugacy $g:S\to\Sigma_3$ is given by 
$g(x)=(...a_{-1},a_0,a_1,...)$ where $a_i=j$ if $f^i(x)\in N_j$.  Notice that
$S\subset\nw(f)$.  Also observe that for points in the  lift, a $0$ in the
itinerary corresponds to movement left and a $2$  corresponds to movement
right.  So, for instance, if $y\in S$ has an  itinerary with a finite number of
$0$s and $1$s and 
$x\in\pi^{-1}(y)$, then $(\tilde f^n(x))_1$ will tend to positive infinity.

Let $y\in S$ be the fixed point with itinerary $(...,0,0,0,...)$ and  let
$x\in\pi^{-1}(y)$.  We claim that $y\in \nw^+(f)\cap\nw^-(f)$ and  moreover,
every sufficiently small disk containing $x$ is both  positively and negatively
returning.  Let $V\subset N$ be the disk 
$\Int(N_0\cap f(N_0))$ and let $U\subset\tilde A$ be the component of 
$\pi^{-1}(V)$ containing $x$.  Examining the dynamics on $\tilde A$  (see Fig.
\ref{fig:triplehorse}) we see that $U$ is negatively  returning (with $n=1$,
$k=-1$) and positively returning (with $n=5$, 
$k=1$).

Moreover, we claim that any disk $W\subset U$ containing $x$ is both 
positively returning and negatively returning.  Let 
$y^\prime=g^{-1}(...,a_0,a_1,...)$ with $a_i=0$ for $i<N$ and $i\ge  3N$ and
$a_i=2$ for $N\le i<3N$.  For $N$ large enough 
$y^\prime,f^{4N-1}(y^\prime)\in\pi(W)$.  Let $x^\prime\in 
W\cap\pi^{-1}(y^\prime)$.  So defined, $\tilde f^{4N-1}(x^\prime)\in  W+1$
(according to the itinerary $x^\prime$ moves left $2N-1$ times  and right $2N$
times).  Thus, $W$ is positively returning with 
$n=4N-1$ and $k=1$.  It is clear that $W$ is negatively returning  with $n=1$,
$k=-1$.
\end{ex}

\section{Fixed points of bounded homeomorphisms}
\label{sec:bddfixed}

In this section we investigate fixed points of bounded homeomorphisms  of the
open annulus.  We begin by applying Theorem 
\ref{thm:connectednw} to this class of homeomorphisms.  We then  describe the
behavior of bounded homeomorphisms possessing one or  fewer fixed points.

\begin{thm}\label{thm:bddfixed} Suppose $f:A\to A$ is a bounded,
orientation-preserving homeomorphism of an open annulus that is homotopic to
the identity, and suppose $\nw(f)$ is connected.  If there is a lift of $f$,
$\tilde f:\tilde A\to\tilde A$, and points $x,y\in\tilde A$ with
$\displaystyle\lim_{n\to\infty}(\tilde f^n(x))_1=-\infty$ and
$\displaystyle\lim_{n\to\infty}(\tilde f^n(y))_1=\infty$, then $\tilde f$, and
hence $f$,has a fixed point.
\end{thm}
\begin{proof} Suppose $x,y\in\tilde A$ with
$\displaystyle\lim_{n\to\infty}(\tilde f^n(x))_1=-\infty$ and
$\displaystyle\lim_{n\to\infty}(\tilde f^n(y))_1=\infty$.  Since $f$ is
bounded, $\omega(\pi(y))\ne\emptyset$.  Let $z\in\omega(\pi(y))$, then let
$y^\prime\in\pi^{-1}(z)$.  If $y^\prime$ is a fixed point then so is $z$, and
we're done.  So assume that $y^\prime$ is not fixed.  Let $U\subset\tilde A$ be
any disk containing $y^\prime$ small enough that $\tilde f(U)\cap U=\emptyset$
and $\pi(U)\subset A$ is a disk. Since $z$ is nonwandering there are infinitely
many positive integers $n$ and corresponding integers $k=k(n)$ such that
$\tilde f^{n}(U)\cap (U+k)\ne\emptyset$.  Since $z\in\omega(\pi(y))$ and
$\displaystyle\lim_{n\to\infty}(\tilde f^n(y))_1=\infty$,  then for $n$ large
enough we can guarantee that $k>0$.  Thus, $U$ is a positively returning disk
with $U\cap\pi^{-1}(\nw(f))\ne\emptyset$. Similarly, since
$\displaystyle\lim_{n\to\infty}(\tilde f^n(x))_1=-\infty$ there is a negatively
returning disk intersecting $\pi^{-1}(\nw(f))$. By Theorem
\ref{thm:connectednw} $f$ has a fixed point.
\end{proof}

In \cite{Ca} Carter considers the case where $g$ is a twist homeomorphism of
the closed annulus $A$ with at most one fixed point in the interior.  She
proves that there is an essential simple closed curve $C$ in the interior which
intersects its image in at most one point.  As we saw in Proposition
\ref{prop:closedannulus}, if $g$ is a bounded homeomorphism of the open
annulus, then there are essential simple closed curves which do not intersect
their images.  Thus it is not clear how one would generalize her theorem for
bounded homeomorphisms.  We do find that bounded homeomorphisms having having
at most one fixed point do have special properties.  We present them in Theorem
\ref{thm:atmostone}.  In particular, we see that if $f$ has at most one fixed
point then the bad behavior found in Example \ref{ex:triplehorse} cannot
occur.  

We state the following theorem for bounded homeomorphisms of the open or closed
annulus.  Recall that for the closed annulus every homeomorphism is bounded;
thus for the closed annulus, the boundedness hypothesis is redundant.

\begin{thm}\label{thm:atmostone} Suppose $f:A\to A$ is an
orientation-preserving, bounded homeomorphism of the open or closed annulus
that is homotopic to the identity, and suppose $f$ has at most one fixed
point.  Let $\tilde f:\tilde A\to\tilde A$ be a lift of $f$.  Then, for each
$x\in\tilde A$ one of the following is true:
\begin{enumerate}
	\item 
	$\displaystyle\lim_{n\to\infty}(\tilde f^n(x))_1=\infty$,
	\item 
	$\displaystyle\lim_{n\to\infty}(\tilde f^n(x))_1=-\infty$, or
	\item 
	$\displaystyle\lim_{n\to\infty}\tilde f^n(x)=p$ for some fixed point $p$ of
$\tilde f$.
\end{enumerate} Moreover, if $\Fix(\tilde f)=\emptyset$ and $\nw(f)$ is
connected, then $\displaystyle\lim_{n\to\infty}(\tilde f^n(x))_1=\infty$ for
all $x\in\tilde A$ or $\displaystyle\lim_{n\to\infty}(\tilde f^n(x))_1=-\infty$
for all $x\in\tilde A$.
\end{thm}
\begin{proof} First, assume that $A$ is the open annulus.  Suppose $f$ has at
most one fixed point.  Let $x\in\tilde A$.  Suppose that $\omega(x)$ is not
empty and consists of more than a fixed point.  Since
$\omega(x)\subset\nw(\tilde f)$, $\nw(\tilde f)$ must consist of more than just
fixed points.  Lemma \ref{lem:nwfp} states that $\tilde f$ has a fixed point of
positive index.  But Corollary \ref{cor:indexzero} states that the Lefschetz 
index of $\Fix(\tilde f)$ is zero; this is a contradiction. Thus
$\omega(x)=\emptyset$ or $\displaystyle\lim_{n\to\infty}\tilde f^n(x)=p$ for
some fixed point $p$ of $\tilde f$.

Now suppose $\omega(x)=\emptyset$.  Since $f$ is bounded Proposition
\ref{prop:closedannulus} states that there is an essential closed annulus
$A_0\subset A$ that is a forward invariant window for $f$.  Let $\tilde
A_0=\pi^{-1}(A_0)$.  Notice that for all $n$ sufficiently large $\tilde
f^n(x)\in\tilde A_0$.  Since we are concerned with the long-term behavior of
$x$, we may assume without loss of generality that $x\in \tilde A_0$.   Since
$\omega(x)=\emptyset$, for any $M>0$ there exists $N_M>0$ such that $|(\tilde
f^n(x)-x)_1|>M$ for all $n>N_M$.  Thus the orbit of $x$ tends to infinity,
negative infinity, or conceivably both.  We will show that the last possibility
will never occur.  Since $A_0$ is compact there is an $M^\prime>0$ such that
$|(\tilde f(y)-y)_1|<2M^\prime$ for all $y\in\tilde A_0$.  Thus, $(\tilde
f^n(x)-x)_1>M^\prime$ for all $n>N_{M^\prime}$ or $(\tilde
f^n(x)-x)_1<-M^\prime$ for all $n>N_{M^\prime}$. So, it must be the case that
$\displaystyle\lim_{n\to\infty}(\tilde f^n(x))_1=\infty$ or
$\displaystyle\lim_{n\to\infty}(\tilde f^n(x))_1=-\infty$.  

Lastly, suppose $\Fix(\tilde f)=\emptyset$ and $\nw(f)$ is connected.  From
above we see that for $x\in\tilde A$ either
$\displaystyle\lim_{n\to\infty}(\tilde f^n(x))_1=\infty$ or
$\displaystyle\lim_{n\to\infty}(\tilde f^n(x))_1=-\infty$.  But, by Theorem
\ref{thm:bddfixed} we know that both cannot occur.

Now, suppose $A$ is the closed annulus.  Then let
$A^\prime=S^1\times(-\ep,1+\ep)$.  Extend $f$ to a bounded homeomorphism on
$A^\prime$ as follows.  If $(x,y)\in S^1\times(1,1+\ep)$, then
$f(x,y)=f(x,1)+(0,(y-1)/2)$.  Similarly define $f$ on $S^1\times(-\ep,0)$. 
Applying the result for the open annulus we arrive at the desired conclusions. 
\end{proof}

\section{Periodic orbits and rotation numbers}
\label{sec:periodic}

As indicated in the introduction, bounded homeomorphisms on noncompact spaces
behave in many ways like homeomorphisms on compact spaces.  In \cite{F3} Franks
proves the following result for homeomorphisms of the closed annulus: if a
point has a given rational rotation number, then there is a periodic point with
that same rotation number.  The result clearly fails for homeomorphisms of the
open annulus.  However, it does hold for bounded homeomorphisms.  

Below we have a theorem that applies to the open and closed annulus.  As
mentioned above, the result for the closed annulus was proved by Franks
(Corollary 2.5 in \cite{F3}) and Handel \cite{H}.   The outline of our proof is
similar to Franks' proof.  However, the results leading up to his proof were
different from those presented here (his arguments used the idea of chain
recurrence), thus we state both results.  In the next two theorems we consider
bounded homeomorphisms of the open and closed annulus.  Recall that for the
closed annulus boundedness is a redundant notion; every homeomorphism of the
closed annulus is bounded.

\begin{thm} Suppose $f:A\to A$ is an orientation-preserving, bounded
homeomorphism of the open or closed annulus that is homotopic to the identity. 
If $\tilde f:\tilde A\to\tilde A$ is a lift of $f$, and for some $x\in\tilde A$
\[\lim\inf\frac{1}{n}(\tilde
f^n(x)-x)_1\le\frac{p}{q}\le\lim\sup\frac{1}{n}(\tilde f^n(x)-x)_1,\] then $f$
has a periodic point with rotation number $p/q$.
\end{thm}
\begin{proof} First, suppose $A$ is the open annulus.  Let $x\in\tilde A$ be a
point satisfying the hypotheses of the theorem.  First, assume that $p=0$.  We
will show that $\tilde f$ has a fixed point.  For the sake of contradiction,
assume that $\tilde f$ has no fixed points.  Then by Theorem
\ref{thm:atmostone} $\lim (\tilde f^n(y)-y)_1=\pm\infty$ for all $y\in \tilde
A$.  Without loss of generality, assume that $\lim(\tilde f^n(x)-x)_1=\infty$. 
Since $f$ is bounded $\omega(\pi(x))\ne\emptyset$.  Let $\tilde\Lambda$ denote
the set $\pi^{-1}(\omega(\pi(x)))$.  So $\lim(\tilde f^n(y)-y)_1=\infty$
uniformly for all $y\in\tilde\Lambda$.  Thus there exists $N>0$ such that
$(\tilde f^N(y)-y)_1>2$ for all $y\in\tilde\Lambda$.  Since the orbit of
$\pi(x)$ limits upon $\pi(\tilde\Lambda)$, for all $k$ sufficiently large
$(\tilde f^{N+k}(x)-\tilde f^k(x))_1>1$. A telescoping sum shows that \[(\tilde
f^{nN+k_0}(x)-\tilde f^{k_0}(x))_1>n\] for some $k_0>0$. Thus,
\[\lim\inf\frac{1}{n}(\tilde f^n(x)-x)_1>\frac{1}{N},\] a contradiction.  Thus
$\tilde f$ has a fixed point.

Now, assume that $p/q\ne 0$.  Let $T:\tilde A\to\tilde A$ be the translation
$T(x,y)=(x+1,y)$.  Let $\tilde g=T^{-p}\circ\tilde f^q$.  So defined, $\tilde
g$ is a lift of $f^q$.  Moreover, $y\in\tilde A$ is a fixed point of $\tilde g$
iff $\pi(y)$ is a periodic point of $f$ with rotation number $p/q$.  Lastly,
observe that \[\lim\inf\frac{1}{n}(\tilde g^{n}(x)-x)_1\le
0\le\lim\sup\frac{1}{n}(\tilde g^{n}(x)-x)_1.\]   Thus, by the argument above
$\tilde g$ has a fixed point, and $f$ has a periodic point with rotation number
$p/q$.

Now, suppose $A$ is the closed annulus.  As in the proof of Theorem
\ref{thm:atmostone} we may extend $f$ to a bounded homeomorphism of the open
annulus $A^\prime=S^1\times(-\ep,1+\ep)$ in such a way that no new periodic
points are created.  Applying the result for the open annulus we find the
prescribed periodic point in $A$. 
\end{proof}

Thus, obviously, if a point has a rational rotation number then there is a
periodic point with the same rotation number.  In fact, we may make the
following conclusion.  The result for the closed annulus was proved by Franks
(\cite{F3})).

\begin{cor} Suppose $f:A\to A$ is an orientation-preserving, bounded
homeomorphism of the open or closed annulus that is homotopic to the identity. 
If among all the periodic points there are only a finite number of periods,
then every point of $A$ has a rotation number.
\end{cor}
\bibliographystyle{cite}
\bibliography{annulus}

\end{document}